\documentclass[
11pt,reqno]{amsart}

\theoremstyle{definition}
\newtheorem{ntn}{Notation}[section]
\newtheorem{dfn}[ntn]{Definition}
\theoremstyle{plain}
\newtheorem{lem}[ntn]{Lemma}
\newtheorem{prp}[ntn]{Proposition}
\newtheorem{thm}[ntn]{Theorem}
\newtheorem{cor}[ntn]{Corollary}
\newtheorem{conj}[ntn]{Conjecture}
\theoremstyle{remark}
\newtheorem{rem}{Remark}

\def\floor[#1]{\lfloor #1 \rfloor }

\newcommand{\N}{\mathbb{N}}
\newcommand{\z}{\mathbb{Z}}
\newcommand{\F}{\mathbb{F}}
\newcommand{\pp}{\mathbb{P}}
\newcommand{\q}{\mathbb{Q}}
\newcommand{\R}{\mathbb{R}}
\newcommand{\C}{\mathbb{C}}

\newcommand{\lan}{\langle}
\newcommand{\ran}{\rangle}

\newcommand{\sr}{{\rm sr}(R)}
\newcommand{\hht}{{\rm ht}}
\newcommand{\rn}{R^{2n}}
\newcommand{\GL}{\mathit{GL}}
\newcommand{\Spp}{\mathit{Sp}}

\newcommand{\iur}{\mathcal{IU}(R^{2n})}
\newcommand{\ivr}{\mathcal{IV}(R^{2n})}

\newcommand{\iurr}{\underline{\mathcal{IU}}(R^{2n})}
\newcommand{\iu}{\mathcal{IU}}

\renewcommand{\o}{\mathcal{O}}
\renewcommand{\H}{\tilde{H}}

\renewcommand{\u}{\mathcal{U}}
\newcommand{\inn}{{\rm inn}}
\newcommand{\inc}{{\rm inc}}

\renewcommand{\ker}{{\rm ker}}

\newcommand{\s}{\Sigma}
\newcommand{\si}{\sigma}

\newcommand{\ph}{\varphi}

\newcommand{\del}{\delta}
\newcommand{\ep}{\epsilon}
\newcommand{\lam}{\Lambda}
\newcommand{\Ga}{\Gamma}

\newcommand{\arr}{\rightarrow}
\newcommand{\larr}{\longrightarrow}

\newcommand{\se}{\subseteq}
\newcommand{\bs}{\backslash}
\newcommand{\bcu}{\bigcup}

\newcommand{\mt}{\mapsto}

\newcommand{\unit}{\mathit{U}^\epsilon_{2n}(R,\Lambda)}

\newcommand{\usr}{{\rm usr}(R)}

\newcommand{\T}{\mathcal{T}}
\renewcommand{\top}{{\rm top}}

\newcommand{\sii}{\overline{\si_p}}
\newcommand{\rr}{{R^\ast}}
\newcommand{\stab}{{\rm Stab}_{G_n}}
\newcommand{\stabe}{{\rm Stab}}
\newcommand{\diag}{{\rm diag}}

\DeclareMathAlphabet{\mathds}{U}{dsrom}{m}{n}
\DeclareMathAlphabet{\mathsc}{U}{rsfs}{m}{n}

\begin{document}

\title{Homology Stability for unitary groups II}
\author{B. Mirzaii}

\begin{abstract}
In this note the homology stability  problem for
hyperbolic unitary groups over a local ring
with an infinite residue field is studied.
\end{abstract}

\maketitle
\maketitle

\section{Introduction}

In this note we continue the study of the homology stability problem for
hyperbolic unitary groups, started in \cite{m-vdk}.
In \cite{m-vdk} a general statement about homology stability for
these groups was established.
It was believed that similar to the general linear group case, for
hyperbolic unitary groups 
over an infinite field, one can have a better range of homology
stability. But no proof of this exists in the literature.
Our main  goal is the study of this problem.


Our main theorem asserts that for a local ring with
an infinite residue field, the natural map
$H_l(\inc): H_l(G_{n}, \z) \arr H_l(G_{n+1}, \z)$
is surjective for $n \ge l+1$ and is injective for $n \ge l+2$, where
$G_n:=\unit$ always with the underlying hyperbolic form.
With a field $k$ as the
coefficient group we get even better result;
$H_l(\inc): H_l(G_{n}, k) \arr H_l(G_{n+1}, k)$
is surjective for $n \ge l$ and is injective for $n \ge l+1$.
In fact the first result follows from
the second one.

To get the second result, we will introduce some posets similar to
one introduced and studied in  \cite{m-vdk}. In section 1 we  prove
that they are highly acyclic.
Applying this we will come up with a spectral
sequence, Theorem \ref{s-vdk}, which is the main purpose
of section 3. The main difficulty is to analyze this
spectral sequence which is done in section 4. The stability theorem
will be a result of this analysis. An application of the stability theorem 
is given in this section.
In section 5, we will discuss the homology stability problem in
the case of a finite field.


I would like to thank W. van der Kallen for his useful
comments that made some of the original proofs shorter 
and for his help in some of the proofs.

Here we will establish some notations. By a ring $R$, we will
always mean a local ring with an infinite residue field
unless it is mentioned. The
ring $R$ has an involution (which may be the identity) and we set
$R_1:=\{r \in R:\overline{r}=r\}$.  This is also a local ring with
an infinite residue field. For the definition of the concepts that
we use such as a bilinear form $h$,
a hyperbolic unitary group and its elementary group,
an isotropic element or set, the unimodular poset $\u(R^n)$, the isotropic
unimodular poset $\iur$ etc, we refer to
\cite[sections 6 and 7]{m-vdk}. We denote a hyperbolic unitary
group $\unit$ and its elementary group by $G_n$ and $E_n$
respectively. By convention $G_0$ will be the trivial group. The
embeddings  $G_n \arr G_{n+1}$ and $E_n \arr E_{n+1}$ are given by
$A \mt \diag(I_2, A)=\left( \begin{array}{cc}
I_2 & 0     \\
0   & A
\end{array} \right)$.
For a group $G$, by $H_i(G)$ we mean $H_i(G, k)$, where $k$ is a field
with trivial $G$-action.
By $k$ as the coefficient group of the homology functor, 
we always mean a field.
In some cases, which will be mentioned, it has to be prime field.

\section{Isotropic unimodular posets}

The main statement of this section, Theorem \ref{iur}, is rather
well known (see \cite{pan1}). We give the details of the proof to
make sure that everything is working for our case, Theorem \ref{iurr}.
For an alternative proof  in the case of a field different from
$\F_2$ see Remark \ref{rem1} and Theorem \ref{finite1}.

\begin{dfn}
Let $S=\{v_1, \dots, v_k\}$ and $T=\{w_1, \dots, w_{k'}\}$ be basis of two
isotropic free summands of $\rn$. We say that $T$ is in general position with
$S$, if $k\le k'$ and the $k' \times k$-matrix $(h(w_i, v_j))$ has a left
inverse.
\end{dfn}

\begin{prp}\label{gp}
Let $n \ge 2$ and assume $T_i$, $1\le i \le l$, are finitely many
finite subsets of $\rn$ such that each $T_i$ is a basis of a free isotropic
summand of $\rn$ with $k$ elements, where $k \le n-1$.
Then there is a basis, $T=\{w_1,\dots, w_n\}$, of a free
isotropic summand of $\rn$ such that $T$ is in
general position with all $T_i$. Moreover ${\rm dim}(W\cap V_i^\perp)=n-k$,
where $W=\lan T \ran$ and $V_i=\lan T_i \ran$.
\end{prp}
\begin{proof}
The proof of the first part is by induction on $l$.
Let $T_i=\{v_{i, 1},\dots, v_{i, k}\}$. For $l=1$,
take a basis of a free isotropic direct summand of $\rn$,
for example $\{w_1, \dots, w_k\}$,
such that $h(w_j, v_{1, m})=\del_{j, m}$, where $\del_{j, m}$
is the Kronecker delta
and choose
$T$ an  extension of this basis to a basis of a maximal isotropic free
subspace. Assume that the
claim is true for $1 \le i \le l-1$. This
means that there is a basis $\{u_1, \dots, u_n\}$ of a free isotropic
summand of $\rn$, in general position with
$T_i$, $1 \le i \le l-1$. Let $\{ x_1, \dots, x_{k}\}$ be a basis of a
free isotropic summand of $\rn$ such that $h( x_j, v_{l, m})=\del_{j, m}$
and take $A=\prod E_{r, s}(a) \in E_n$
such that $Au_j=x_j$, $1 \le j \le k$ \cite[6.5, 7.1]{m-vdk}.
If $x_j:=Au_j$ for $k+1 \le j \le n$, then
$\{ x_1, \dots, x_n\}$ is in general position with
$T_l$. Set $B_i=(h(u_j, v_{i, m}))$, $1 \le i \le l-1$ and
$B_l=(h(x_j, v_{l,m}))$. Let $B_i^{(k)}$ be the matrix obtained from
$B_i$ by deleting $(j_{1, i}, \dots, j_{n-k, i})$-th rows such that
$f_i^{(k)}:={\rm det}(B_i^{(k)}) \in R^\ast$ for all $i$.
Set $A(t)=\prod E_{r, s}(ta)$, $B_i(t)=(h(A(t)u_j, v_{i, m}))$,
for $1 \le i \le l$ and let $B_i^{(k)}(t)$ be the matrix obtained
by deleting $(j_{1, i}, \dots, j_{n-k, i})$-th rows
of $B_i(t)$ and set $f_i^{(k)}(t):={\rm det}(B_i^{(k)}(t)) \in R[t]$.
Clearly $f_i^{(k)}(0)=f_i^{(k)}$ for $1 \le i \le l-1$ and
$f_l^{(k)}(1)=f_l^{(k)}$. It is not difficult to see that there is a
$t_1 \in R$ such that $f_i^{(k)}(t_1) \in R^\ast$ for $1 \le i \le l$
\cite[1.4, 1.5]{kal2}. Take $W=\{A(t_1)u_1, \dots, A(t_1)u_n \}$.
The second part of the proposition follows from the exact sequence
\[
0 \arr W\cap V_1^\perp \arr W \overset{\psi}{\arr} R^k \arr 0
\]
with $\psi(w):=(h(w, v_{1, 1}), \dots, h(w, v_{1, k}))$ and
the fact that projective modules over local rings are free.
\end{proof}

Let $S$ be a non-empty set and $X \se \o(S)$ \cite[Sec. 4]{m-vdk}. Let
$C_k(X)$, $k\ge 0$, be the free
$\z$-module with the basis consisting of the $k$-simplices ($(k+1)$-frames)
of $X$, $C_{-1}(X)=\z$ and $C_k(X)=0$ for $k \le -2$. The family
$C_\ast(X):=\{C_k(X)\}$ yields a chain complex with the differentials
$\partial_0: C_0(X) \arr C_{-1}(X)=\z$, $\sum_i n_iv_i \mt \sum_i n_i$ and
$\partial_k=\sum_{i=0}^k(-1)^id_i: C_k(X) \arr C_{k-1}(X)$, $k\ge 1$, where
$d_i((v_0, \dots, v_k))= (v_0, \dots,\widehat{v_i}, \dots, v_k)$.
The poset $X$ is called $n$-acyclic if $\H_k(X,\z):= H_k(C_\ast(X))=0$
for all $0\le k \le n$.

\begin{lem}\label{tttt}
Let $n, m$ be two natural numbers and $n \le m$. If $n \ge k+1$ then for
every $(v_1, \dots, v_k) \in \u(R^m)$ there is a
$v\in R^n=\lan e_1, \dots, e_n \ran$
such that $(v, v_1, \dots, v_k) \in \u(R^m)$.
\end{lem}
\begin{proof}
The proof is similar to the proof of lemma \cite[5.4]{m-vdk}, using the
fact that $\sr=1$.
\end{proof}

\begin{thm}\label{ur}
Let $n, m$ be two natural numbers and $n \le m$. Then the poset
$\o(R^n)\cap \u(R^m)$ is $(n-2)$-acyclic and $\o(R^n)\cap \u(R^m)_w$ is
$(n-|w|-2)$-acyclic for every $w=(w_1, \dots, w_r) \in \u(R^m)$.
\end{thm}
\begin{proof}
Let $X=\o(R^n) \cap \u(R^m)=\u(R^n)$ and
$\si=\sum_{i=1}^l n_i(v_0^i, \dots, v_k^i)$ be a cycle
in $C_{k}(X)$, $k \le n-2$. It is not difficult to see that there is a
unimodular vector $v \in R^n$ such that
$\{v, v_0^i, \dots, v_k^i\}$ is linearly independent, $1 \le i \le l$. If
$\beta:=\sum_{i=1}^l n_i(v, v_0^i, \dots, v_k^i) \in C_{k+1}(X)$,
then $\partial_{k+1}(\beta)=\si$, so $X$ is $(n-2)$-acyclic.

Let $Y=\o(R^n)\cap \u(R^m)_w$ and assume that $n-|w|-2\ge -1$.
Let $\si$ be a $k$-cycle in $C_k(Y)$
with $k \le n-|w|-2$.
To prove the second part of the theorem  it is sufficient
to find a  unimodular vector $v \in R^n$ such that
$\{v, v_0^i, \dots, v_k^i, w_1, \dots, w_r\}$ is
linearly independent, $1 \le i \le l$. The proof
is by induction on $l$. The case $l=1$ follows from \ref{tttt}.
By induction assume that there are
$u_1, u_2 \in R^n$ such that
$(u_1, v_0^i, \dots, v_k^i, w_1, \dots, w_r )\in \u(R^m)$
for $1 \le i \le l-1$ and
$(u_2, v_0^l, \dots, v_k^l, w_1, \dots, w_r) \in \u(R^m)$.
Let $A=\prod E_{r, s}(a)$ be an element of the elementary group
$E_n(R) \se \GL_n(R)$ such that
$Au_1=u_2$ and set $A(t)=\prod E_{r, s}(ta)$. Let $B_i$ be the matrix whose
columns are the vectors $u_1, v_0^i, \dots, v_k^i, w_1, \dots, w_r$ for
$1 \le i \le l-1$, $B_l$ the matrix whose columns are
$u_2, v_0^l, \dots, v_k^l, w_1, \dots, w_r$ and $B_i(t)$ is the
matrix whose columns are $A(t)u_1, v_0^i, \dots, v_k^i, w_1, \dots, w_r$,
$1 \le i \le l$. The rest of the proof is similar to the proof of
proposition \ref{gp}.
\end{proof}

\begin{thm}\label{iur}
The poset $\iur$ is $(n-2)$-acyclic.
\end{thm}
\begin{proof}
If $n=1$, then everything is trivial, so we assume that $n \ge 2$.
Let $\si=\sum_{i=1}^rn_iv_i$ be a $k$-cycle.
Thus $v_i$, $1 \le i \le r$, are isotropic
$(k+1)$-frames with $k \le n-2$. 
By \ref{gp},
there is an isotropic $n$-frame $w$ in general position with 
$v_i$, $1 \le i \le r$. Set $W=\lan w \ran$ and let $E_\si$ be the set of all
$(u_1, \dots, u_m, t_1, \dots, t_l) \in \iur$ such that $m, l \ge 0$,
$(u_1, \dots, u_m) \in \u(W)$, if $m \ge 1$,
and for every $l \ge 1$ there exist an $i$
such that  $(t_1, \dots, t_l) \le v_i$.
The poset $E_\si$ satisfies the chain condition and
$v_i \in E_\si$.
It is sufficient to prove that $E_\si$ is $(n-2)$-acyclic, because then
$\si \in \partial_{k+1}(E_\si) \se \partial_{k+1}(C_{k+1}(X))$.
Let $F:= E_\si$. Since
$\o(W) \cap F=\u(W)$, by \ref{ur} the poset  $\o(W) \cap F$ is $(n-2)$-acyclic.
If $u \in F\bs \o(W)$, then $u$ is of the
form $(u_1, \dots, u_m, t_1, \dots, t_l)$, $l\ge 1$. By \ref{gp},
${\rm dim}(V)=n-l$, where
$V=W \cap \lan t_1 \dots t_l\ran^\perp$. With all this we have
\begin{gather*}
\o(W)\cap F_u=\o(V)\cap \iur_{(u_1, \dots, u_m)}=
\o(V)\cap \u(W)_{(u_1, \dots, u_m)}.
\end{gather*}
Again by \ref{ur}, $\o(V)\cap \u(W)_{(u_1, \dots, u_m)}$ is
$((n-l)-m-2)$-acyclic, so $\o(W)\cap F_u$ is
$(n-|u|-2)$-acyclic. Therefore $F$ is $(n-2)$-acyclic \cite[2.13 (i)]{kall}.
\end{proof}

\begin{rem}\label{rem1}
\par (i) The concept of {\it being in general position} and the idea
of the proof of \ref{iur} is taken from \cite{pan1}.
Because the details of the proof in \cite{pan1} never appeared
we wrote it down.
\par (ii) In fact Theorem \ref{iur} is true for every
field $R \neq \z/2\z$. Let
\[
\ivr=\{ V \se \rn: V \neq 0 {\rm \ and \ isotropic \ subspace}\}.
\]
Define the map of the posets $f:\iur \arr \ivr$, $v\mt \lan v \ran$.
As Vogtmann proved, \cite[Thm. 1.6]{vog}, $\ivr$ is $(n-2)$-connected
(Vogtmann proved this for $G_n=O_{2n}(R)$, but her proof works without
modification in our more general setting \cite[p. 115]{ch}).
On the one hand it is easy to see that
$Link_{\ivr}^+(V) \simeq \mathcal{IV}(R^{2(n-{\rm dim}(V))})$, so it is
$(n-{\rm dim}(V)-2)$-connected and on the other hand
$f/V=\u(V)$ which is $({\rm dim}(V)-2)$-connected \cite[2.6]{kall},
hence defining the height function
$\hht_{\ivr}(V)={\rm dim}(V)-1$ \cite[section 2]{m-vdk},
one sees that $\iur$ is $(n-2)$-connected
\cite[Thm. 3.8]{m-vdk}.
\par (iii) We expect that over a ring with no finite ring as a
homomorphic image and finite
unitary stable rank the poset $\iur$ is $(n-\usr -1)$-connected.
For this it is sufficient to prove \ref{gp} over such ring. For example
Theorem \ref{iur}, without any change in its proof, is true over a semi-local
ring with infinite residue
fields. Therefore the results of this note are also valid for these rings.
\par (iii) Using \ref{iur}, (iii)  and the same argument as in (ii)
one can prove that over a semi-local
ring with infinite residue fields, $\ivr$ is $(n-2)$-acyclic.
Over an infinite field this gives  much easier proof of Vogtmann's 
theorem mentioned in (ii).
\par (iv) Using a theorem of Van der Kallen \cite[Thm. 2.6]{kall}
and a similar arguments as (iii) we can generalize the Tits-Solomon
theorem over a ring with stable range one (for example any Artinian ring).
Let $R$ be a ring with stable range one  and consider the following poset,
which we call it the Tits poset,
\[
\T(R^n)=\{ V \se R^n: V  \ {\rm free \ summand \ of} \ R^n , V \neq
0, R^n\}.
\]
Let $X=\u(R^n)_{\le n-1}$ and consider the poset map $g: X \arr \T(R^n)$,
$v \mt \lan v\ran$. By induction and a similar argument as in (ii),
using the fact that $X$ is $(n-3)$connected,
one can prove that $\T(R^n)$ is $(n-3)$-connected (note that
any stably free projective module of rank $\ge 1$ is free).
We leave the details of the proof to the interested readers.
\end{rem}

\begin{dfn}
Define $\underline{\u}(R^n)=\{ (\lan v_1 \ran, \dots, \lan v_k \ran):
(v_1, \dots, v_k) \in \u(R^n)\}$ and
$\iurr=\{ (\lan v_1 \ran, \dots, \lan v_k \ran): (v_1, \dots, v_k) \in
\iur\}$.
\end{dfn}

\begin{thm}\label{iurr}
Let $n, m$ be two natural numbers and $n \le m$.
Then the poset $\o(\pp^{n-1})\cap \underline{\u}(R^m)$ is $(n-2)$-acyclic,
the poset $\o(\pp^{n-1})\cap \underline{\u}(R^m)_w$ is $(n-|w|-2)$-acyclic
for every $w \in \underline{\u}(R^m)$ and  the  poset $\iurr$ is
$(n-2)$-acyclic.
\end{thm}
\begin{proof}
The proof is similar to the proof of \ref{ur} and \ref{iur}.
\end{proof}

\section{The spectral sequence}

In this section, $k$ will be a field, $S_i$
a $k$-algebra, $i \in \N$, $S_i^{\otimes n}:= S_i \otimes_k \dots \otimes_k S_i$
($n$-times) and $V_{n}(S_i):=(S_i^{\otimes n})^{\s_n}$, where
$\s_n$ is the symmetric group of degree $n$.

\begin{lem}\label{n-s1}
Let $\ph_i: R \arr S_i$ be a ring homomorphism, $i=1, \dots, d$.
Consider the action of $R^\ast$ on
$\bigotimes_{i=1}^d S_i^{\otimes n_i}$
and $\bigotimes_{i=1}^d V_{n_i}(S_i)$ as
\begin{gather*}
r\bigotimes_{i=1}^d(a_{1, i}\otimes \dots \otimes a_{n_i, i})=
\bigotimes_{i=1}^d(\ph_i(r)^{t_i}a_{1, i}\otimes \dots \otimes
\ph_i(r)^{t_i}a_{n_i, i}),
\end{gather*}
where $t_i\ge 1$.
Then  $H_0(R^\ast,\bigotimes_{i=1}^d S_i^{\otimes n_i})=
H_0(R^\ast, \bigotimes_{i=1}^d V_{n_i}(S_i))=0$.
\end{lem}
\begin{proof}
The proof is similar to the proof of \cite[1.5]{nes-sus}
and \cite[1.6]{nes-sus} with minor generalization.
If $B:= \bigotimes_{i=1}^d S_i^{\otimes n_i}$, then
$H_0(\rr, B)=B/I$, where $I$ is the ideal of $B$ generated by the elements
$\bigotimes_{i=1}^d(\ph_i(r)^{t_i}\otimes \dots \otimes \ph_i(r)^{t_i})-1$.
Consider the collection $\{ \psi_1^{(j_i)}, \dots, \psi_{t_i}^{(j_i)} \}$,
$i=1,\dots, d$, $1 \le j_1 \le n_1$ and
$\sum_{i=1}^{m-1}n_i < j_m \le \sum_{i=1}^{m}n_i $ for $m \ge 2$,
of homomorphisms $R \arr B/I$ given by
$\psi_l^{(j_i)}(r)=1\otimes \dots \otimes \ph_i(r) \otimes
\dots \otimes 1\mod I$,
with $\ph_i(r)$  in the $j_i$-th position.
For simplicity we denote this collection by
$\{ \psi_l: 1 \le l \le \sum_{i=1}^{d}t_in_i
\}$. If $I$ is a proper ideal, we
obtain a collection of ring homomorphisms $\psi_l$
such that $\prod_{l}\psi_l(r)= 1$
for every $r \in \rr$, but this is
impossible \cite[Cor. 1.3, Lem. 1.4]{nes-sus}.
Thus $I=B$ and therefore $H_0(R^\ast,\bigotimes_{i=1}^d S_i^{\otimes n_i})=0$.
For the proof of the second part let
$l_1^{(i)}, \dots, l_{s_i}^{(i)}$, $i=1, \dots, d$, be
the natural numbers such that $\sum_{j=1}^{s_i} l_j^{(i)}=n_i$, and
denote by $V_{n_i}^{l_1^{(i)}, \dots, l_{s_i}^{(i)}}$ the subspace of
$V_{n_i}(S_i)$ generated by the elements of the form
\begin{gather*}
y_{c, l^{(i)},i}:=
\sum_{\del \in \s_{n_i}/\s_{l_1^{(i)}}
\times \dots \times \s_{l_{s_i}^{(i)}}}
(\underset{l_1^{(i)}}
{\underbrace{c_1^{(i)}\otimes \dots \otimes c_1^{(i)}}}
\otimes \dots \otimes
\underset{l_{s_i}^{(i)}}{\underbrace{c_{s_i}^{(i)}\otimes \dots \otimes
c_{s_i}^{(i)}}})^\del.
\end{gather*}
Clearly  $V_{n_i}^{l_1^{(i)}, \dots, l_{s_i}^{(i)}}$
is an $\rr$-invariant subspace of $V_{n_i}(S_i)$ and
$V_{n_i}(S_i)=
{\sum_{l_1^{(i)}+\cdots +l_{s_i}^{(i)}=n_i}}
V_{n_i}^{l_1^{(i)}, \dots, l_{s_i}^{(i)}}$. Let
$V_{n_i}^{(j)}(S_i)=
\sum_{s_i\ge n_i-j}V_{n_i}^{l_1^{(i)}, \dots, l_{s_i}^{(i)}}$
and set
\[
T_h:=\sum_{h_1 + \cdots +h_d=h}
V_{n_1}^{(h_1)}(S_1)\otimes \cdots
\otimes V_{n_d}^{(h_d)}(S_d).
\]
It is not difficult to
see that if $\sum_{i=1}^{d}n_i-s_i=h$ and
$l_1^{(i)}+\cdots +l_{s_i}^{(i)}=n_i$, then
\[
\bigotimes_{i=1}^d S_i^{\otimes s_i}
\arr T_h/T_{h-1} ,
\]
\[
\bigotimes_{i=1}^d c_1^{(i)} \otimes \dots \otimes c_{s_i}^{(i)}
\mt y_{c, l^{(1)},1} \otimes \cdots
\otimes y_{c, l^{(d)},d} \mod T_{h-1}
\]
is multilinear, so it gives an $\rr$-equivariant homomorphism. In this way we
obtain an $\rr$-equivariant epimorphism
\[
\coprod_{n_1-s_1+ \cdots +n_d-s_d=h
}\quad
\bigotimes_{i=1}^d S_i^{\otimes s_i} \quad\arr\quad T_h/T_{h-1}.
\]
Since the functor $H_0$ is right exact,  by applying  the first part of the
lemma we get $H_0(\rr, T_h/T_{h-1})=0$. By induction on $h$ we prove that
$H_0(\rr, T_h)=0$. If $h=0$, then
$T_0=\bigotimes_{i=1}^d V_{n_i}^{(0)}(S_i)$ and
$\bigotimes_{i=1}^d S_i^{\otimes n_i} \arr T_0$  is surjective, so
$H_0(\rr, T_0)=0$. By induction and applying the functor $H_0$ to the short exact
sequence $0 \arr T_{h-1} \arr T_h \arr T_h/T_{h-1} \arr 0$, we see that
$H_0(\rr, T_h)=0$.
\end{proof}

\begin{lem}\label{n-s2}
Let $P_i$ and $Q_i$ be $S_i$-modules for $i=1, \dots, d$. Then
$\bigotimes_{i=1}^d \bigwedge^{n_{i, 1}}P_i \otimes_k V_{n_{i, 2}}(Q_i)$
has a natural structure of $\bigotimes_{i=1}^d V_{n_i}(S_i)$-module,
where $n_i=n_{i, 1}+n_{i, 2}$. Moreover for all $l \ge 0$
\begin{gather*}
H_l(\rr, \bigotimes_{i=1}^d \bigwedge\nolimits^{n_{i, 1}}P_i \otimes_k
V_{n_{i, 2}}
(Q_i))=0.
\end{gather*}
\end{lem}
\begin{proof}
The first part follows immediately from \cite[Lem. 1.7]{nes-sus} and
the second part follows from \ref{n-s1} and \cite[Lem. 1.8]{nes-sus}.
\end{proof}

Let $B$ be a $k$-vector space and let $\Ga(B)$ be
the algebra of divided powers of
$B$, which is a graded commutative algebra concentrated
in even degrees and endowed with a system of divided powers
with $\Ga_{2n}(B)=V_n(B)$ (see (\cite[Chap. V, No. 6]{bro} and
\cite[\S 1]{nes-sus} for more details). The homology of an abelian
group $A$ with rational coefficients coincides
with exterior powers: $H_p(A, \q)= \bigwedge^{p}(A\otimes \q)$. The
homology with coefficients in the prime field
$\F_p=\z/p\z$ is more complicated. The ring $H_\ast(A, \F_p)$
has a canonical structure of divided powers
\cite[Chap. V, Example 6.5.4]{bro}.
Moreover, $H_1(A, \F_p)=A/pA$ and there is an exact sequence
\[
0 \arr \bigwedge\nolimits^2 (A/pA) \arr H_2(A, \F_p) \overset{\beta}{\arr}
{}_pA \arr 0.
\]
Any choice of a section for $\beta$ gives a homomorphism
$\ph :{}_pA \arr H_2(A, \F_p)$, which by the property of the
algebra $\Ga$, uniquely extends to an $\F_p$-algebra homomorphism
$\bigwedge\nolimits^2 (A/pA) \otimes_k \Ga({}_pA) \arr H_\ast(A, \F_p)$,
thus giving rise to
an isomorphism of graded $\F_p$-algebras \cite[Chap. V, Thm. 6.6]{bro},
\cite[ \S 8, Prop. 3]{quil6}.
We identify $H_j(A, \F_p)$ with
$\coprod_i \bigwedge^{j-2i} (A/pA) \otimes_k \Ga_{2i}({}_pA)$ and introduce
a filtration on $H_j(A, \F_p)$ setting
\[
H_j^{(r)}= \coprod_{i \le r}\bigwedge\nolimits^{j-2i} (A/pA)
\otimes_k \Ga_{2i}({}_pA).
\]
This filtration does not depend on our choice of section $\ph$
and successive factors $H_j^{(r)}/H_j^{(r-1)}$ are canonically isomorphic
to $\bigwedge^{j-2r} (A/pA) \otimes_k \Ga_{2r}({}_pA)$.

\begin{thm}\label{n-s3}
Let $M_i$ be a $T_i$-module and let $R \arr T_i$ be a ring homomorphism.
Consider the action of $R^\ast$ on $M_i$ given by
$r\cdot m=\ph_i(r)^{t_i}m$, where $t_i\geq1$. If $k$ is a prime field,
then $H_l(R^\ast, \bigotimes_{i=1}^d H_{l_i}(M_i))=0$
for  $l\ge 0$, where $l_i > 0$ for some $i$.
\end{thm}
\begin{proof}
Let $P_i=M_i\otimes_\z k$ and $S_i=T_i\otimes_\z k$. If $k=\q$, then
$H_{l_i}(M_i)= \bigwedge^{l_i} P_i$
and if $k=\F_p$ for some prime number $p$,
then we can find an $\rr$-invariant filtration  of
$\bigotimes_{i=1}^d H_{l_i}(M_i)$ whose successive
factors are isomorphic to
$\bigotimes_{i=1}^d \bigwedge^{j_i-2m_i}P_i \otimes_k \Ga_{2m_i}(Q_i)$
for some $j_i$ and $m_i$, where $Q_i={}_p(P_i\otimes_\z k)$. Note that
$P_i$ and $Q_i$ are $S_i$-modules. Then both cases follow from \ref{n-s2}.
\end{proof}

Let $\overline{\si_2}=(\lan e_1 \ran, \lan e_3 \ran) \in \iurr$.
The elements of
$\stabe_{G_n}(\overline{\si_2})=\{ B \in G_n: B\overline{\si_2}=
\overline{\si_2}\}$ are of the form
\[
\left(\begin{array}{cccccc}
a_1 &  \ast   &  0   &  \ast    &  \ast  &  \ast         \\
 0  & {\overline{a_1}}^{-1} &  0   &    0     &   0    &    0          \\
 0  & \ast     & a_2  &  \ast    &  \ast  &  \ast         \\
 0  &  0       &  0   & {\overline{a_2}}^{-1} &   0    &    0          \\
 0  &  \ast    &  0   &  \ast    &        &               \\
 0  &  \ast    &  0   &  \ast    &        &    A
\end{array} \right),
\]
where $a_i \in \rr$ and $A \in G_{n-2}$. Let $N_{n, 2}$ and
$L_{n, 2}$ be  the subgroups of $\stabe_{G_n}(\overline{\si_2})$ of elements
of the form
\begin{center}
$\left(\begin{array}{cccccc}
 1  &  \ast   &  0   &  \ast    &  \ast  &  \ast         \\
 0  &   1      &  0   &    0     &   0    &    0          \\
 0  & \ast     &  1   &  \ast    &  \ast  &  \ast         \\
 0  &  0       &  0   &    1     &   0    &    0          \\
 0  &  \ast    &  0   &  \ast    &        &               \\
 0  &  \ast    &  0   &  \ast    &        &  I_{2(n-2)}
\end{array} \right)$,
$\left(\begin{array}{cccccc}
a_1 &  \ast   &  0   &  \ast    &  \ast  &  \ast          \\
 0  & {\overline{a_1}}^{-1} &  0   &    0     &   0    &    0          \\
 0  & \ast     & a_2  &  \ast    &  \ast  &  \ast         \\
 0  &  0       &  0   & {\overline{a_2}}^{-1} &   0    &    0          \\
 0  &  \ast    &  0   &  \ast    &        &               \\
 0  &  \ast    &  0   &  \ast    &        &   I_{2(n-2)}
\end{array} \right)$
\end{center}
respectively. It is a matter of an easy calculation
to see that the elements of the group
${N}_{n, 2}'=[N_{n, 2},N_{n, 2}]$ are of the form
\[
\left(\begin{array}{cccccc}
 1  &   r      &  0   &     t    &   0    &    0         \\
 0  &   1      &  0   &     0    &   0    &    0         \\
 0  &   -\ep^{-1}\overline{t}   &  1   &     s    &   0    &    0    \\
 0  &   0      &  0   &     1    &   0    &    0         \\
 0  &   0      &  0   &     0    &        &              \\
 0  &   0      &  0   &     0    &        &   I_{2(n-2)}
\end{array} \right),
\]
where $r,\ s \in \lam=\{r \in R: \ep^{-1}\overline{r}=-r\}$ and $t \in R$.
In general one can define $N_{n, p}$,
$L_{n, p}$ and $N_{n, p}'$ for all $p$, $1 \le p \le n$, in a similar way.
Embed  $\rr^p \times G_{n-p}$ in $\stab(\sii)$ as
$ \diag(a_1, \dots, a_p, A) \mt
\diag(\left(
\begin{array}{cc}
a_1 & 0          \\
0   & {\overline{a_1}}^{-1}
\end{array}
\right)
, \dots,
\left(
\begin{array}{cc}
a_p & 0          \\
0   & {\overline{a_p}}^{-1}
\end{array}
\right), A)$.

\begin{thm}\label{m-vdk5}
Let $\sii=(\lan e_1 \ran, \lan e_3 \ran,
\dots, \lan e_{2p-1}\ran) \in \iurr$.
Then the inclusion $\rr^p \times G_{n-p} \arr \stab(\sii)$ induces
the isomorphism between the homology groups
$H_i(\rr^p \times G_{n-p}) \arr H_i(\stab(\sii))$ for all $i$.
\end{thm}
\begin{proof}
It is sufficient to prove the theorem when $k$ is a prime field.
Fix a natural number $p$, $1 \le p \le n$, and set $N=N_{n, p}$, $L=L_{n, p}$,
$N'=N_{n, p}'$ and $T=\stab(\sii)$.
The extensions $1 \arr N' \arr L \arr L/N' \arr 1$ and
$1 \arr N/N' \arr L/N' \arr L/N \arr 1$ give 
the Lyndon-Hochschild-Serre spectral sequences
\begin{gather*}
\hspace{-4.3 cm}
E_{p, q}^2=H_p(L/N', H_q(N')) \Rightarrow H_{p+q}(L),\\
E_{p', q'}^2=H_{p'}(L/N, H_{q'}(N/N', H_q(N'))) \Rightarrow
H_{p'+ q'}(L/N', H_q(N')),
\end{gather*}
respectively. Since $L/N\simeq \rr^p$ and $N/N'$
acts trivially on $N'$, $E_{p', q'}^2=H_{p'}(\rr^p,
H_{q'}(N/N')\otimes_k H_q(N'))$.  
It is not difficult to see that $N/N'\simeq R^h$ and
$N'\simeq R^l \times \lam^m$ for some $h, l, m$ and the action of
$R_1^*$ on $N/N'$ and $N'$ is linear-diagonal and
quadratic-diagonal respectively. Again 
the extension $1 \arr
R_1^\ast  \arr \rr^p \arr \rr^p/R_1^\ast \arr 1$ ($R_1^\ast$
embeds in $\rr^p$ diagonally) gives
\[
E_{r, s}^2=H_r(\rr^p/R_1^\ast, H_s(R_1^\ast, M)) \Rightarrow
H_{r+s}(\rr^p, M),
\]
where $M= H_{q'}(N/N')\otimes_k H_q(N')$.
Since the homology functor commutes with the direct sum functor,
\[
H_s(R_1^\ast, M)\simeq \bigoplus_{i=0}^q H_s(R_1^\ast,
H_{q'}(R^h) \otimes_k H_{i}(R^l) \otimes_k H_{q-i}(\lam^m)),
\]
where the action of $R_1^*$ on $R^h$, $R^l$ and $\lam^m$  is linear-diagonal,
quadratic-diagonal and quadratic-diagonal respectively.
By theorem \ref{n-s3},
$H_s(R_1^\ast, M)=0$  for  $s \ge 0$ and $q>0$ or $q'>0$.
This shows that $E_{p', q'}^2=0$ for  $p'\ge 0$ and $q>0$ or $q'>0$.
Therefore $H_{p'}(L/N', H_q(N'))=0$ for $p' \ge 0$ and $q>0$.
Hence $E_{p, q}^2=0$ for $p \ge 1$ and $q >0$.
By the convergence of the spectral sequence we get
\begin{equation}\label{isom1}
H_p(L) \overset{\simeq}{\arr} H_{p}(L/N').
\end{equation} 
The extension $1 \arr N/N' \arr L/N' \arr L/N \arr 1$ gives
\[
E_{i, j}^2=H_{i}(L/N, H_{j}(N/N')) \Rightarrow
H_{i+j}(L/N').
\]
and by a similar approach to (\ref{isom1}),
\begin{equation}\label{isom2}
H_i(\rr^p) \overset{\simeq}{\arr} H_{i}(L/N').
\end{equation}
From the embedding $\rr^p \arr L$,
(\ref{isom1}) and (\ref{isom2}) we get the
 isomorphism $H_i(\rr^p) \overset{\simeq}
 {\arr} H_{i}(L)$, $i\ge 0$.
The commutative diagram
\[
\begin{array}{ccccccccc}
 1  & \arr & \rr^p     & \arr  & \rr^p \times G_{n-p} &   \arr   &
 G_{n-p}     & \arr &  1  \\
    &      &\downarrow &       & \downarrow           &          &
    \downarrow  &      &     \\
 1  & \arr &  L        & \arr  &  T                   &   \arr   &
 G_{n-p}     & \arr &  1
\end{array}
\]
gives the map of the spectral sequences
\[
\begin{array}{ccc}
E_{p, q}^2 = H_p(G_{n-p}, H_q(\rr^p)) & \Rightarrow &
H_{p+q}(\rr^p \times G_{n-p})  \\
\Big\downarrow                        &             &  \Big\downarrow  \\
\!\!\!\!\!
{E'}_{p, q}^2 = H_p(G_{n-p}, H_q(L))  & \Rightarrow & H_{p+q}(T).
\end{array}
\]
By what we proved in the above we have the isomorphism
$ E_{p, q}^2 \simeq  {E'}_{p, q}^2$.
This gives an isomorphism on the abutments and so $
H_{i}(\rr^p \times G_{n-p}) \simeq H_{i}(T)$.
\end{proof}

\begin{thm}\label{s-vdk}
There is a first quadrant spectral sequence converging to zero with
\[
E_{p, q}^1(n)= \begin{cases} H_q(\rr^p \times G_{n-p}) &
\text{if $0 \le p \le n$}\\
H_q(G_n, H_{n-1}(X_n)) & \text{if $p=n+1$} \\
0 & \text{if $p \ge n+2$}\end{cases},
\]
where $X_n=\iurr$.\\ For $1 \le p \le n$ the differential
$d_{p, q}^1(n)$ equals
$\sum_{i=1}^p(-1)^{i+1}H_q(\alpha_{i, p})$, where
$\alpha_{i, p}: \rr^p \times G_{n-p} \arr \rr^{p-1} \times G_{n-p+1}$,
$ \diag(a_1, \dots, a_p, A) \mt
\diag(a_1, \dots, \widehat{a_i}, \dots, a_p,
\left(
\begin{array}{ccc}
a_i & 0         & 0  \\
0   & {\overline{a_i}}^{-1}  & 0  \\
0   &  0        & A
\end{array}
\right))$.\\ In particular for $0 \le p \le n$,
$d_{p, 0}^1(n)= \begin{cases} {\rm id}_k &
\text{if $p$ is odd}\\
0 & \text{if $p$ is even}\end{cases}$,
so $E_{p, 0}^2=0$ for $0 \le p \le n-1$.
\end{thm}
\begin{proof}
Let $C_l'(X_n)$ be the $k$-vector space
with the basis consisting of $l$-simplices (isotropic $(l+1)$-frames)
of $X_n$. Since $X_n$ is $(n-2)$-acyclic, Theorem \ref{iurr},
we get an exact sequence
\[
0 \leftarrow k \leftarrow C_0'(X_n) \leftarrow C_1'(X_n) \leftarrow \cdots
\leftarrow C_{n-1}'(X_n) \leftarrow  H_{n-1}(X_n, k) \leftarrow 0.
\]
Call this exact sequence  $L_\ast$:
$L_0= k, L_i=C_{i-1}'(X_n)$, $1 \le i \le n$,
$L_{n+1}= H_{n-1}(X_n, k)$ and $L_i=0$ for $i\ge n+2$.
Let $F_\ast \arr k$ be a resolution of $k$ by free (left) $G_n$-modules
and consider the bicomplex $C_{\ast,\ast}=L_\ast \otimes_{G_n} F_\ast$.
Here we convert the left action of $G_n$ on $L_\ast$ into a right action
with $vg:=g^{-1}v$. By the
general theory of the spectral sequence for a bicomplex we have
$E_{p, q}^1(I)=H_q(C_{p, \ast})=H_q(L_p \otimes_{G_n} F_\ast)$ and
$E_{p, q}^1(\mathit{II})=H_q(C_{\ast, p})=H_q(L_\ast \otimes_{G_n} F_p)$. 
Since $F_p$ is a free $G_n$-module, $L_\ast \otimes_{G_n} F_p$ is exact
and this shows that $E_{p, q}^1(\mathit{II})=0$. 
Therefore $E_{p, q}^1(n):= E_{p, q}^1(I)$ converges to zero.
If $p=0$, then $E_{0, q}^1(n)=H_q(k \otimes_{G_n} F_\ast)=H_q(G_n)$.
The group $G_n$ acts transitively
on the $l$-frames of $X_n$, $1 \le l \le n$,
so by the Shapiro lemma \cite[Chap. III, 6.2]{bro},
$L_p \otimes_{G_n} F_\ast \simeq k \otimes_{\stab(\sii)} F_\ast$ and thus
$E_{p, q}^1(n)=H_q(\stab(\sii))$, $1 \le p \le n$.
By \ref{m-vdk5},
$E_{p, q}^1(n)=H_q(\rr^p \times G_{n-p})$ for $0 \le p \le n$,
hence $E_{p, q}^1(n)$ is of the form that we mentioned.
Now we look at the differential
$d_{p, q}^1(n): E_{p, q}^1(n) \arr E_{p-1, q}^1(n)$, $1 \le p \le n$;
$d_{1, q}^1(n)$ is induced by 
$C_0'(X_n)\otimes_{G_n} F_\ast \arr k \otimes_{G_n} F_\ast$,
$\overline{\si_1} \otimes x \mt  1\otimes x$.
Considering the isomorphism
$k \otimes_{\stabe_{G_n}(\overline{\si_1})} F_\ast \arr C_0'(X)
\otimes_{G_n} F_\ast$, 
$1 \otimes x \mt  \overline{\si_1}\otimes x$, one sees that $d_{1, q}^1(n)$
is induced by
$k \otimes_{\stabe_{G_n}(\overline{\si_1})}
F_\ast \arr k \otimes_{G_n} F_\ast$.
This shows that $d_{1, q}^1(n)$
is the map $H_q(\stabe_{G_n}(\overline{\si_1})) \arr H_q(G_{n})$
induced by the inclusion map, therefore
$d_{1, q}^1(n)=H_q(\inc):
H_q(\inc): H_q(\rr \times G_{n-1}) \arr H_q(G_{n})$.
For $2 \le p \le n$, $d_{p, q}^1(n)$
is induced by the
map $\sum_{i=1}^{p}(-1)^{i+1}d_i: L_{p} \arr L_{p-1}$,
where  $d_i$ deletes the
$i$-th component of the isotropic $p$-frames. Let $g_{i, p}$ be the permutation
matrix such that
$(e_{2h-1}, e_{2h})g_{i, p}^{-1}=(e_{2h-1}, e_{2h})$, $1 \le h \le i-1$,
$(e_{2i-1}, e_{2i})g_{i, p}^{-1}=(e_{2p-1}, e_{2p})$ and
$(e_{2l-1}, e_{2l})g_{i, p}^{-1}=(e_{2l-3}, e_{2l-2})$, $i+1 \le l \le p$,
where $vg^{-1}:=gv$ for $v \in \rn$.
It is easy to see that $d_i(\sii)=\overline{\si_{p-1}}g_{i, p}$, and
so $\partial(\sii)=\sum_{i=1}^p(-1)^{i+1}\overline{\si_{p-1}}g_{i, p}$.
Consider 
$d_i\otimes {\rm id}_{F_\ast}:
L_p \otimes_{G_n} F_\ast \arr L_{p-1} \otimes_{G_n} F_\ast$,
$\sii \otimes x \mt d_i(\sii) \otimes x$.
Let $\inn_{g_{i, p}}:G_n \arr G_n$, $g \mt g_{i, p}gg_{i, p}^{-1}$ and
$l_{g_{i, p}}: F_\ast \arr F_\ast$, $x \mt g_{i, p}x$.
It is easy to see that $l_{g_{i, p}}$ is an
$\inn_{g_{i, p}}$-homomorphism, and $d_i\otimes {\rm id}_{F_\ast}$
induces the map
$k \otimes_{\stab(\sii)} F_\ast \arr
k \otimes_{\stabe_{G_n}(\overline{\si_{p-1}})} F_\ast$,
$1 \otimes x \mt   1\otimes l_{g_{i, p}}(x)$.
This shows that $d_i$ induces
$H_q(\inn_{g_{i, p}}):H_q(\stab(\sii)) \arr
H_q(\stabe_{G_n}(\overline{\si_{p-1}}))$
and hence the map
$H_q(\inn_{g_{i, p}}): H_q(\rr^p \times G_{n-p}) \arr
H_q(\rr^{p-1} \times G_{n-p+1})$.
Set $\alpha_{i, p}= \inn_{g_{i, p}}$.
Since  $G_n$ acts transitively on
the generators of ${C_p'}(X_n)$,
$E_{\ast,0}^1(n)$
is of the following form
\[
0 \leftarrow k \leftarrow k \leftarrow k \leftarrow \cdots
\leftarrow k \leftarrow H_0(G_n H_{n-1}(X_n, k)) \leftarrow 0,
\]
where
$d_{p, 0}^1(n)= \begin{cases} {\rm id}_k &
\text{if $p$ is odd}\\
0 & \text{if $p$ is even}\end{cases}.$
Clearly $E_{p, 0}^2(n)=0$ if $0 \le p \le n-1$.
\end{proof}

\begin{rem}
In fact $E_{n, 0}^2(n)=0$. For a proof see the proof
of theorem \ref{m3}.
\end{rem}

\section{Stability theorem}

To prove the homology stability result
we have to study the spectral
sequence that we obtained in theorem \ref{s-vdk}.

\begin{lem}\label{m1}
Let $n \ge 1$, $l \ge 0$ be integer numbers such that $n-1 \ge l$.
Let  $H_q(\inc): H_q(G_{n-2}) \arr H_q(G_{n-1})$
be surjective if $0 \le q \le l-1$. Then the following conditions are
equivalent;
\par {\rm (i)}  $H_l(\inc): H_l(G_{n-1}) \arr H_l(G_{n})$
is surjective,
\par {\rm (ii)} $H_l(\inc): H_l(\rr \times G_{n-1})
\arr H_l(G_{n})$ is surjective.
\end{lem}
\begin{proof}
For $n=1$ every thing is easy
so let $n \ge 2$.
By the K\"unneth theorem \cite[Chap. V, \S 10, Thm. 10.1]{mac}
we have 
$ H_l(\rr \times G_{n-1})= S_1 \oplus S_2$, where
$S_1= H_l(G_{n-1})$ and
$S_2= \bigoplus_{i=1}^l H_i(\rr)\otimes_k H_{l-i}(G_{n-1})$.
The case (i)$\Rightarrow$(ii) is trivial. To prove
(ii)$\Rightarrow$(i) it is sufficient to prove that
$\tau_1(S_2) \se \tau_1(S_1)$, where $\tau_1=H_l(\inc)$.
{}From $ i \ge 1$ and $n-1 \ge l$, we have
$n-2 \ge l-1 \ge l-i$, so by hypothesis
$H_{l-i}(\inc): H_{l-i}(G_{n-2}) \arr H_{l-i}(G_{n-1})$ is surjective,
$1 \le i \le l$. Consider the following diagram
\[
\begin{array}{ccccc}
H_i(\rr) \otimes_k H_{l-i}(G_{n-1}) & \overset{\beta_1}{\arr} &
H_l(\rr \times G_{n-1}) & \overset{\tau_1}{\arr} & H_l(G_{n})\\
\Big\uparrow \vcenter{%
\rlap{$\scriptstyle{\alpha_1}$}}
              &              &  &
&\Big\uparrow\vcenter{%
\rlap{$\scriptstyle{\alpha_2}$}}  \\
H_i(\rr) \otimes_k H_{l-i}(G_{n-2}) &\overset{\beta_2}{\arr} &
H_l(\rr \times G_{n-2}) & \overset{\tau_1'}{\arr} & H_l(G_{n-1}),
\end{array}
\]
where $\beta_j$ is the shuffle product, $j=1,2$ \cite[Chap. V, Sec. 5]{bro},
$\alpha_1={\rm id}\otimes H_{l-i}(\inc)$ is surjective
and $\alpha_2=H_l(\inc)$. By giving an explicit description of the
above maps we prove that this diagram is commutative.
For this purpose we use
the the bar resolution of a group \cite[Chap. I, Sec. 5]{bro}. If
$x=\sum [a_1| \dots| a_i]\otimes [A_1| \dots| A_{l-i}]
\in H_i(\rr) \otimes_k H_{l-i}(G_{n-2})$, then
\begin{gather*}
\hspace{-3.5 cm}
x \overset{\alpha_1}{\longmapsto}
\sum [a_1| \dots| a_i]\otimes [\diag(I_2, A_1)|
\dots| \diag(I_2, A_{l-i})] \\
\overset{\tau_1\circ \beta_1}{\longmapsto}
\sum\sum_\del
{{\rm sign}(\del) }
[\dots| \diag(a_{\del(i')}, {\overline{a_{\del(i')}}}^{-1},
I_{2(n-1)})| \dots| 
\diag(
I_4, A_{\del(j')})| \dots]
\end{gather*}
and
\begin{gather*}
\hspace{-0.9 cm}
x \overset{\tau_1' \circ \beta_2}{\longmapsto}
\sum\sum_\del
{{\rm sign}(\del)}[\dots|
\diag(a_{\del(i')},{\overline{a_{\del(i')}}}^{-1},
I_{2(n-2)})| \dots |\\
\hspace{8 cm}
 \diag(I_2, A_{\del(j')})| \dots] \\
\overset{\alpha_2}{\longmapsto}
\sum\sum_\del
{{\rm sign}(\del) }
[\dots| \diag(I_2, a_{\del(i')},
{\overline{a_{\del(i')}}}^{-1}, I_{2(n-2)})|\dots |\\
\hspace{7.9 cm}
 \diag(
 I_4, A_{\del(j')})| \dots].
\end{gather*}
See \cite[Chap. VIII, \S 8]{mac} for more details about the shuffle product.
Let $P \in G_n$ be the permutation matrix that permutes the
first and second columns
with third and forth columns respectively and let $\inn_P: G_n \arr G_n$,
$A \mt PAP^{-1}=PAP$. It is well known that
$H_q(\inn_P)={\rm id}_{H_q(G_n)}$ \cite[Chap. II, \S 8]{bro}, hence
\begin{gather*}
H_l(\inn_P)(
[\dots| \diag(I_2, a_{\del(i')},
{\overline{a_{\del(i')}}}^{-1}, I_{2(n-2)})|
\dots| \diag(I_2, I_2, A_{\del(j')})| \dots]) \\
\ \ \ \ \ \ \ \ \ \
=[\dots| \diag(a_{\del(i')},
{\overline{a_{\del(i')}}}^{-1}, I_2, I_{2(n-2)})|
\dots| \diag(I_2, I_2, A_{\del(j')})| \dots].
\end{gather*}
This shows that the above diagram is commutative. Therefore
$\tau_1(S_2) \se \tau_1(S_1)$.
\end{proof}

\begin{lem}\label{m2}
Let $n \ge 2$, $l \ge 0$ be integer numbers such that $n-1 \ge l+1$.
Let $H_q(\inc ): H_q(G_{m-1}) \arr H_q(G_{m})$ be isomorphism
for $m=n-1, n-2$ and $0 \le q \le {\min}\{l-1, m-2\}$.
Then the following conditions are equivalent;
\par {\rm (i)} $H_l(\inc): H_l(G_{n-1}) \arr H_l(G_{n})$ is bijective,
\par {\rm (ii)} $H_l(\rr^2\times G_{n-2}) \overset{\tau_2}{\arr}
H_l(\rr \times G_{n-1}) \overset{\tau_1}{\arr} H_l(G_{n}) \arr 0$ is exact,
where $\tau_1=H_l(\inc)$ and $\tau_2= H_l(\alpha_{1,2})- H_l(\alpha_{2,2})$.
\end{lem}
\begin{proof}
Let $H_l(\rr \times G_{n-1})= S_1 \oplus S_2$, where
$S_1$ and $S_2$ are as in the proof of lemma \ref{m1} and
$H_l(\rr^2 \times G_{n-2})=\bigoplus_{h=1}^4 T_h$, where
\begin{gather*}
\hspace{-7.2 cm}
T_1= H_l(G_{n-2}), \\
\hspace{-3.7 cm}
T_2= \bigoplus_{i=1}^l H_i(\rr \times 1)\otimes_k H_{l-i}(G_{n-2}),\\
\hspace{-3.7 cm}
T_3= \bigoplus_{i=1}^l H_i(1\times \rr)\otimes_k H_{l-i}(G_{n-2}), \\
\!\!\!\!\!\!\!\!\!\!
T_4= \bigoplus_
{i+j \le l}
H_i(\rr \times 1)\otimes_k H_j(1\times \rr)\otimes_k
H_{l-i-j}(G_{n-2}).
\end{gather*}
Set $\si_1^{(2)}=H_l(\alpha_{1,2})$ and $\si_2^{(2)}=H_l(\alpha_{2,2})$.
First  (i)$\Rightarrow$(ii).
The surjectivity
of $\tau_1$ is trivial. Let $(x , v) \in S_1 \oplus S_2$ such that
$\tau_1((x , v))=0$. The relations $n-1 \ge l+1$ and $i \ge 1$ imply
that $n-3 \ge l-1 \ge l-i$ and hence
$H_{l-i}(G_{n-2}) \arr H_{l-i}(G_{n-1})$ is bijective,
so there exists $w \in T_2$ such that $-\si_1^{(2)}(w)=v$. If
$y=(0, w, 0, 0) \in \bigoplus_{h=1}^4 T_h$, then
$\tau_2(y)=(\si_2^{(2)}(w), -\si_1^{(2)}(w))=(\si_2^{(2)}(w), v)$.
Since $\tau_1\circ \tau_2=0$,
$\tau_1(\si_2^{(2)}(w))=-\tau_1(v)$. Combining this with
$\tau_1(x)=-\tau_1(v)$, we obtain $\tau_1(\si_2^{(2)}(w))=\tau_1(x)$.
By injectivity of $H_l(\inc)$,  $\si_2^{(2)}(w)=x$, thus
$\tau_2(y)=(x, v)$. This shows that
the complex is exact. The proof of (ii)$\Rightarrow$(i) is more difficult.
The surjectivity of $\tau_1=H_l(\inc)$ follows from lemma \ref{m1}.
Let $ x \in {\rm ker}(H_l(G_{n-1}) \arr H_l(G_{n}))$,
then $(x, 0) \in {\rm ker}(\tau_1)$. By exactness of the complex there is a
$y=(0, y_2, y_3, y_4) \in \bigoplus_{h=1}^4 T_h$ such that $\tau_2(y)=(x, 0)$
(one should notice that $\tau_2(T_1)=0$).
First we prove that we can assume $y_4=0$. If $n=2$, then $l \le1$ and so
$T_4=0$. Therefore we may assume $n \ge 3$.
Consider the summand
\begin{gather*}
U=%
{\bigoplus_{{i,j\ge 1}}}
H_i(\rr \times 1 \times 1)\otimes_k H_j(1 \times \rr\times 1)\otimes_k
H_{l-i-j}(G_{n-3})
\end{gather*}
of $H_l(\rr^3 \times G_{n-3} )$ and set
$\tau_3:=d_{3, l}^1(n)=\si_1^{(3)}-\si_2^{(3)}+\si_3^{(3)}$,
where $\si_i^{(3)}=H_l(\alpha_{i, 3})$.
It is easy to see that $\si_3^{(3)}(U) \se T_4$ and
$-\si_2^{(3)}+\si_1^{(3)}(U) \se T_2 \oplus T_3$. By assumption,
$\si_3^{(3)}|_U: U \arr T_4$ is an isomorphism. If $\si_3^{(3)}(z)=y_4$,
then $y-\tau_3(z)=(0, y_2', y_3', 0)$
and $\tau_2(y-\tau_3(z))=(x, 0)$.
So we can assume that $y=(0, y_2, y_3, 0)$. Let
\begin{gather*}
y_2=(\sum [a_1| \dots| a_i]\otimes [A_1| \dots| A_{l-i}])_{1 \le i \le l}\\
y_3=(\sum [b_1| \dots| b_i]\otimes [B_1| \dots| B_{l-i}])_{1 \le i \le l}.
\end{gather*}
By an explicit computation
\[
\tau_2(y)=(\si_1^{(2)}(y_2)-\si_2^{(2)}(y_3),
-\si_2^{(2)}(y_2)+\si_1^{(2)}(y_3)).
\]
This shows that $x=\si_1^{(2)}(y_2)-\si_2^{(2)}(y_3)$
is equal to
\begin{gather*}
\sum_{i=1}^l\sum
\sum_{\del}
{{\rm sign}(\del)}[\dots| \diag(a_{\del(i')},
{\overline{a_{\del(i')}}}^{-1}, I_{2(n-1)})| \dots| 
\diag(I_2, A_{\del(j')})|
\dots]\\
-\sum_{i=1}^l\sum
\sum_{\del}
{{\rm sign}(\del)}[\dots| \diag(b_{\del(i')},
{\overline{b_{\del(i')}}}^{-1}, I_{2(n-1)})| \dots| 
\diag(I_2, B_{\del(j')})|\dots]
\end{gather*}
and for $1 \le i \le l$,
\begin{gather*}
0=\sum [a_1| \dots| a_i]\otimes [\diag(I_2, A_1)|
\dots| \diag(I_2, A_{l-i})] \\
\ \ \ \
-\sum [b_1| \dots| b_i]\otimes [\diag(I_2, B_1)|
\dots| \diag(I_2, B_{l-i})].
\end{gather*}
By the injectivity of $H_{l-i}(G_{n-2}) \arr H_{l-i}(G_{n-1})$,
we see that $y_2=y_3$, (note that we view $y_2$ and $y_3$ as elements
of  $T_2$ or $T_3$). Now it is easy to see that $x=0$.
\end{proof}

Consider $R^{2(n-2)}$ as the submodule of $\rn$ generated
by $e_5, e_6, \dots, e_{2n}$ (so $G_{n-2}$ embeds in $G_n$ as
$\diag(I_2, I_2, G_{n-2})$). Let
${L'}_\ast$ be the complex
\begin{gather*}
\hspace{5 cm}
\cdots \leftarrow C_{n-3}'(X_{n-2}) \leftarrow  H_{n-3}(X_{n-2},k)
\leftarrow 0 \\
\hspace{-4 cm}
0\leftarrow 0\leftarrow 0 \leftarrow k \leftarrow C_0'(X_{n-2})
\leftarrow C_1'(X_{n-2})
\leftarrow
\end{gather*}
with $X_{n-2}=\underline{\iu}(R^{2(n-2)})$.
Define the map of complexes
${L'}_\ast \overset{\alpha_\ast}{\arr} {L}_\ast$, given by
\begin{gather*}
(\lan v_1\ran, \dots, \lan v_k\ran) \overset{\alpha_k}{\mt}
(\lan e_1\ran,\lan e_3\ran,\lan v_1\ran, \dots, \lan v_k\ran)-
(\lan e_1\ran, \lan e_1+e_3\ran,\lan v_1\ran, \dots, \lan v_k\ran)\\
\!\!\!\!\!\!\!\!\!\!\!\!\!\!
+(\lan e_3\ran, \lan e_1+e_3\ran, \lan v_1\ran, \dots, \lan v_k\ran)
\end{gather*}
Note that this is similar to one defined in the proof of the proposition
2.6 in \cite{nes-sus}. This gives the maps of  bicomplexes
\begin{gather*}
{L'}_\ast \otimes_{G_{n-2}} {F'}_\ast \arr
{L}_\ast \otimes_{G_{n}} {F}_\ast \arr
{L}_\ast \otimes_{G_{n}} {F}_\ast /{L'}_\ast \otimes_{G_{n-2}} {F'}_\ast,
\end{gather*}
where $L_\ast$ and $F_\ast$ are as in the proof of theorem \ref{s-vdk}
and ${F'}_\ast$ is ${F}_\ast$ as $G_{n-2}$-module, so it induces the maps
of spectral sequences
\begin{gather*}
{E'}_{p, q}^r(n) \arr {E}_{p, q}^r(n) \arr {E''}_{p, q}^r(n),
\end{gather*}
where all the three spectral sequences converge to zero.
By a similar argument  as in the proof of \ref{s-vdk},
one sees that the spectral sequence ${E'}_{p, q}^1(n)$ is of the form
\begin{gather*}
{E'}_{p, q}^1(n)= \begin{cases}
E_{p-2, q}^1(n-2) &
\text{if $p \ge 2$}\\
0 & \text{if $p=0, 1$} \end{cases}.
\end{gather*}
For $2 \le p \le n$,
${E'}_{p, q}^1(n) \arr {E}_{p, q}^1(n)$ is 
induced by
$\inc:\rr^{p-2} \times G_{n-p} \arr \rr^{p} \times G_{n-p}$,
$A \mt \diag(I_2, I_2 , A)$, and 
\begin{gather*}
{E''}_{p, q}^1(n)={E}_{p, q}^1(n)/{E'}_{p, q}^1(n).
\end{gather*}
{}From the complexes
\begin{gather*}
\!\!\!\!\!\!\!\!\!\!\!\!
{D}_\ast(q) : \ \ \ \
0 \arr {E}_{n, q}^1(n) \arr {E}_{n-1, q}^1(n)
\arr \cdots \arr {E}_{0, q}^1(n) \arr 0 \\
\!\!\!\!\!\!\!\!
{D'}_\ast(q) : \ \ \ 0 \arr {E'}_{n, q}^1(n) \arr {E'}_{n-1, q}^1(n)
\arr \cdots \arr {E'}_{0, q}^1(n) \arr 0 \\
\!\!\!\!
{D''}_\ast(q) : \ \ \ \! 0 \arr {E''}_{n, q}^1(n) \arr {E''}_{n-1, q}^1(n)
\arr \cdots \arr {E''}_{0, q}^1(n) \arr 0
\end{gather*}
we obtain a short exact sequence
\[
0 \arr {D'}_\ast(q) \arr {D}_\ast(q) \arr {D''}_\ast(q) \arr 0
\]
and by applying  the homology long exact sequence to this short exact
sequence we get the following exact sequence
\begin{gather*}
\hspace{-1.5 cm}
 {E'}_{n-1, q}^2(n)
\arr {E}_{n-1, q}^2(n)
\arr {E''}_{n-1, q}^2(n)
\arr {E'}_{n-2, q}^2(n) \\
\hspace{4 cm}
\arr \cdots
\arr {E'}_{0, q}^2(n)
\arr {E}_{0, q}^2(n)
\arr {E''}_{0, q}^2(n)
\arr
0.
\end{gather*}

\begin{thm}\label{m3}
Let $n \ge 1$, $l \ge 0$ be integer numbers. Then
$H_l(\inc): H_l(G_{n-1}) \arr H_l(G_{n})$
is surjective for $n-1 \ge l$ and is injective for $n-1 \ge l+1$.
\end{thm}
\begin{proof}
The proof is by induction on $l$. If $l=0$ then everything is obvious.
Assume the induction hypothesis, that is $H_i(G_{m-1}) \arr H_i(G_{m})$
is surjective if $m-1 \ge i$ and is bijective if $m-1 \ge i+1$,
where $1 \le i \le l-1$.
Let $n-1 \ge l$ and
consider the spectral sequence ${E''}_{p, q}^2(n)$.
To prove the surjectivity, it is sufficient  to prove that
${E''}_{p, q}^2(n)=0$ if $n \ge p+q$, $0 \le q \le l-1\le n-2$
and $2 \le p \le n-1$,
because then we obtain ${E''}_{0, l}^2(n)=0$ and applying lemma
\ref{m1} we have the desired result. Let $R_i^\ast$ denotes the $i$-th factor
of ${\rr}^m$. By the K\"unneth theorem
${E''}_{p, q}^1(n)=T_1 \oplus T_2 \oplus T_3 \ \mod \ {E'}_{p, q}^1$,
where
\begin{gather*}
T_1= \bigoplus_{i_1\ge 1}
H_{i_1}(R_1^\ast)\otimes H_{i_3}(R_3^\ast) \otimes \cdots \otimes
H_{i_{p}}(R_p^\ast) \otimes H_{q-\s i_j}(G_{n-p}) ,\\
T_2= \bigoplus_{i_2\ge 1}
H_{i_2}(R_2^\ast)\otimes H_{i_3}(R_3^\ast) \otimes \cdots \otimes
H_{i_{p}}(R_p^\ast) \otimes H_{q-\s i_j}(G_{n-p}) ,\\
\ \ \ \ \
T_3= \bigoplus_{k_1,k_2\ge 1}
H_{k_1}(R_1^\ast)\otimes H_{k_2}(R_2^\ast) \otimes \cdots \otimes
H_{k_{p}}(R_p^\ast) \otimes H_{q-\s k_s}(G_{n-p}).
\end{gather*}
Consider the following summand of ${E''}_{p+1, q}^1(n)$
\begin{gather*}
U_1=\bigoplus_{j_2,j_3\ge 1}
H_{j_2}(R_2^\ast)\otimes H_{j_3}(R_3^\ast) \otimes \cdots \otimes
H_{j_{p+1}}(R_{p+1}^\ast) \otimes H_{q-\s j_t}(G_{n-p-1}),
\end{gather*}
where $j_t=k_{t-1}$, $2 \le t \le p+1$.
Let $\si_i^{(m)}:=H_l(\alpha_{i, m})$.
It is easy to see that
$\si_1^{(p+1)}(U_1) \se T_3$ and
$\sum_{i=2}^{p+1}(-1)^{i+1}\si_i^{(p+1)}(U_1) \se T_2$.
Let $x=(x_1, x_2, x_3) \in \ker({d''}_{p, q}^1)$.
Since $n-p-1 \ge q-1\ge q-\sum j_t$,
by a similar argument as in the proof of lemma \ref{m2},
we can assume that $x_3=0$. If
\begin{gather*}
U_2=\bigoplus_{j_2\ge 1}
H_{j_2}(R_2^\ast)\otimes H_{j_4}(R_4^\ast) \otimes \cdots \otimes
H_{j_{p+1}}(R_{p+1}^\ast)\otimes H_{q-\s j_t}(G_{n-p-1}),
\end{gather*}
then we have $\si_1^{(p+1)}(U_2) \se T_1$,
$\si_2^{(p+1)}(U_2)=0 \ {\rm mod} \ {E'}_{p, q}^1(n)$ and
$\sum_{i=3}^{p+1}(-1)^{i+1}\si_i^{(p+1)}(U_2) \se T_2$.
In the same way, using our assumption we can again assume that
$x_1=0$. So $x=(0, x_2, 0)$. Once again we have
$\si_1^{(p)}(T_2) \se S_1$ and
$\sum_{i=2}^{p}(-1)^{i+1}\si_i^{(p)}(T_2) \se S_2$, where
\begin{gather*}
\!\!\!\!
S_1=\bigoplus_{k_1\ge 1}
H_{k_1}(R_1^\ast)\otimes H_{k_2}(R_2^\ast) \otimes \cdots \otimes
H_{k_{p-1}}(R_{p-1}^\ast) \otimes H_{q-\s k_t}(G_{n-p+1}),\\
\!\!\!\!\!\!\!\!\!
S_2=\bigoplus_{l_2\ge 1}
H_{l_2}(R_2^\ast)\otimes H_{l_3}(R_3^\ast) \otimes \cdots \otimes
H_{l_{p-1}}(R_{p-1}^\ast) \otimes H_{q-\s l_t}(G_{n-p+1}).
\end{gather*}
By induction hypothesis $\si_1^{(p)}$ is an isomorphism, so
$x_2=0$. Therefore ${E''}_{p, q}^2(n)=0$
if $n \ge p+q$, $2 \le p \le n-1$, $1 \le q \le l-1$.
To prove that ${E''}_{p, 0}^2(n)=0$ for $0 \le p \le n$, it is sufficient to
prove that ${E}_{p, 0}^2(n)=0$ for $0 \le p \le n$.
For $0 \le p \le n-1$ this follows from \ref{s-vdk}.
If $n$ is odd then  $E_{n, 0}^2(n)=0$, because
$d_{n, 0}^1(n)= {\rm id}_k$. So let $n$ be even. We prove
by induction on $n$ that $E_{n, 0}^2(n)=0$. If $n=2$, then
\[
\theta:= (\lan e_1\ran, \lan e_3\ran) - (\lan e_1\ran, \lan e_1+e_3\ran)
+(\lan e_3\ran, \lan e_1+ e_3\ran) \in H_1(X_2, k)
\]
and so $d_{3, 0}^1(2)(\theta  \mod G_2)= 1 \in \z$. 
Assume that this is true for $n-2$, that is $E_{n-2, 0}^2(n-2)=0$.
From the map 
${E'}_{p, q}^1(n) \arr {E}_{p, q}^1(n)$ we get the commutative
diagram
\begin{gather*}
\begin{array}{ccccc}
H_0(G_{n-2}, H_{n-3}(X_{n-2}, k))  & \overset{d_{n-1, 0}^1(n-2)}
{-\!\!\!-\!\!\!-\!\!\!-\!\!\!-\!\!\!-\!\!\!-\!\!\!\rightarrow}
& k&\longrightarrow & 0\\
\Big\downarrow\vcenter{%
\rlap{$\scriptstyle{\alpha'}$}} &      &
\Big\downarrow\vcenter{%
\rlap{$\scriptstyle{{\rm id}_k}$}} & & \\
H_0(G_{n}, H_{n-1}(X_n, k))  & \overset{d_{n+1, 0}^1(n)}
{-\!\!\!-\!\!\!-\!\!\!-\!\!\!-\!\!\!-\!\!\!-\!\!\!\rightarrow}  & k
& \longrightarrow & 0,
\end{array}
\end{gather*}
where the map $\alpha'$ is induced by the map $\alpha_\ast$.
By induction
and the commutativity of the above diagram we see that $d_{n+1, 0}^1(n)$
is surjective and therefore $E_{n, 0}^2(n)=0$.
This shows that
${E''}_{p, 0}^2(n)=0$, $0 \le p \le n$ and so the proof of 
the claim is complete. To complete the proof of the theorem we must prove the 
injectivity claimed in the theorem. This
can be done by a similar argument as in the above with suitable changes.
%
\end{proof}

\begin{cor}\label{m4}
If $n-p \ge q$, then the complex
\begin{gather*}
\hspace{-2 cm}
H_q(\rr^p \times G_{n-p}) \overset{\tau_p}{\larr}
H_q(\rr^{p-1} \times G_{n-p+1}) \overset{\tau_{p-1}}{\larr} \cdots \\
\hspace{5 cm}
\overset{\tau_2}{\larr} H_q(\rr \times G_{n-1}) \overset{\tau_1}{\larr}
H_q( G_{n}) \larr 0
\end{gather*}
is exact, where $\tau_i:=d_{i, q}^1(n)$.
\end{cor}
\begin{proof}
This comes out of the proof of \ref{m3}.
\end{proof}

\begin{thm}\label{m5}
Let $n \ge 1$, $l \ge 0$ be integer numbers. Then
$H_l(\inc): H_l(G_{n}, \z) \arr H_l(G_{n+1}, \z)$
is surjective for $n \ge l+1$ and is injective for $n \ge l+2$.
\end{thm}
\begin{proof}
For $n \ge l+1$, theorem \ref{m3} implies
$H_{l+1}(G_{n+1}, G_{n})=0$. Here $H_{l+1}(G_{n+1}, G_{n})$ is
the homology of the mapping cone of the map of complexes
$F_\ast^{(n)} \arr F_\ast^{(n+1)}$ with coefficients in k where
$F_\ast^{(m)}$ is the $G_m$-resolution of $k$.
Applying the homology long exact sequence to the short exact sequence
\[
0 \arr \z \arr \q \arr \q/\z \arr 0
\]
we have the exact sequence
\begin{gather*}
\hspace{-3 cm}
\cdots \arr H_{l+1}(G_{n+1}, G_{n}, \q/\z) \arr H_l(G_{n+1}, G_{n}, \z) \\
\hspace{3 cm}
\arr H_l(G_{n+1}, G_{n}, \q) \arr H_l(G_{n+1}, G_{n}, \q/\z) \arr \cdots .
\end{gather*}
We must prove that $H_{l+1}(G_{n+1}, G_{n}, \q/\z)=0$.
Since $\q/\z=\oplus_p\;\underset{\longrightarrow}{{\rm lim}} \z/p^d\z$
and since the homology functor commutes with the
direct limit functor, it is sufficient
to prove that  $H_{l+1}(G_{n+1}, G_{n}, \z/p^d\z)=0$. This can be 
deduced from writing the homology long exact sequence 
of the short exact sequence
\[
0 \arr \z/p\z \arr \z/p^d\z \arr \z/p^{d-1}\z \arr 0
\]
and induction on $d$. Therefore $H_l(G_{n+1}, G_{n}, \z)=0$.
The surjectivity, claimed in the theorem, follows
from the long exact sequence
\begin{gather*}
\hspace{-3.5 cm}
\cdots \arr H_{l+1}(G_{n+1}, G_{n}, \z) \arr H_l(G_{n}, \z) \\
\hspace{3.5 cm}
\arr H_l(G_{n+1}, \z) \arr H_l(G_{n+1}, G_{n}, \z) \arr \cdots.
\end{gather*}
The proof of the other claim follows from a similar argument.
\end{proof}

\begin{rem}
Theorem \ref{m5} gives almost a positive answer to a question
asked by Sah in \cite[4.9]{sah}. Also it gives better range of stability
in comparison to other results \cite{m-vdk}, \cite{vog}.
\end{rem}

Let $G$ be a topological group and let $BG^\top$ be the quotient space
$
\bcu_n \Delta^n \times G^n/ \sim$, where $\sim$ is the relation
\begin{gather*}
\hspace{-8 cm}
(t_0, \dots, t_n,g_1, \dots, g_n) \sim \\
\ \ \ \ \ \ \ \ \ \ \ \
\begin{cases}
(t_0, \dots, \widehat{t_i}, \dots, t_n,
g_1, \dots, g_{i-1}, g_i g_{i+1}, g_{i+2},
\dots, g_n) & \text{if $t_i=0$}\\
(t_0, \dots, t_{i-1}, t_i + t_{i+1}, t_{i+2}, \dots, t_n,
g_1, \dots ,\widehat{g_i},\dots, g_n) & \text{if $g_i=e$} .
\end{cases}
\end{gather*}
It is easy to see that $B$ is a functor from the the category of topological
groups to the category of topological spaces. 
The topological space $BG^\top$ is called the
classifying space of $G$ with the underlying topology.
Let $BG$ be the classifying space of $G$ as the topological group
with the discrete topology. By the
functorial property of $B$ we have a natural map $\psi: BG \arr BG^\top$.

\begin{conj}[Friedlander-Milnor Conjecture]
Let $G$ be a Lie group. The canonical map
$\psi: BG \arr BG^\top$ induces isomorphism of homology and cohomology
with any finite abelian coefficient group.
\end{conj}

See \cite{mil} and \cite{sah} for more information in this direction.

\begin{thm}\label{sus-kar}
Let $F=\R$ or $\C$. If  $G=O(F)$, $\Spp(F)$ or  $U(\C)$, then
$H_i(BG, A) \simeq H_i(BG^\top, A)$ for all $i$ and any finite
coefficient group A.
\end{thm}
\begin{proof}
See \cite[Thm. 1, 2]{kar1}
\end{proof}

\begin{cor}
Let $F=\R$ or $\C$. If  $G_n=O_{2n}(F)$, $\Spp_{2n}(F)$ or $U_{2n}(\C)$,
then $H_i(BG_n, A) \simeq H_i(BG_n^\top, A)$ if $n \ge i+1$ and any finite
coefficient group A.
\end{cor}
\begin{proof}
This follows from \ref{m3} and \ref{sus-kar}.
\end{proof}

\section{Homology stability of unitary groups over finite fields}

In this section we will  explain which part of the above results is
true  if $R$ is a finite field, so in this
section we assume that $R:=F$ is a finite field.

\begin{lem}\label{finite1}
Let $F$ be a field different from $\F_2$. Then $\underline{\u}(F^n)$
is $(n-2)$-connected, $\underline{\u}(F^n)_w$ is $(n-|w|-2)$-connected
for every $w \in \underline{\u}(F^n)$ and  the  poset
$\underline{\iu}(F^{2n})$ is $(n-2)$-connected.
\end{lem}
\begin{proof}
The proof of  the first two claims is by induction
on $n$. Let $Z:= \underline{\u}(F^n)$
and $Y:=\o(\pp^{n-2})$. For any
$v=(\lan v_1 \ran , \dots, \lan v_k \ran) \in Z \bs Y$,
there is an $i$, for example $i=1$, such that $v_i \notin R^{n-1}$.
This means that the $n$-th coordinate of $v_1$ is not zero.
Choose $r_i\in F$ such that $v_i'=v_i-r_iv_1 \in F^{n-1}$, $2 \le i \le k$.
It is not difficult to see that
\begin{gather*}
Y \cap Z_v\simeq Y \cap
\underline{\u}(F^n)_{(\lan v_1 \ran ,\lan v_2' \ran, \dots, \lan v_k' \ran)}
\simeq Y \cap
\underline{\u}(F^n)_{(\lan v_2' \ran, \dots, \lan v_k' \ran)} \\
\hspace{-3.57 cm}
\simeq \underline{\u}(F^{n-1})_{(\lan v_2' \ran, \dots, \lan v_k' \ran)}.
\end{gather*}
By induction
$\underline{\u}(F^{n-1})_{(\lan v_2' \ran, \dots, \lan v_k' \ran)}$
is $((n-1)-(|v|-1)-2)$-connected, so $Y \cap Z_v$ is
$((n-3)-|v|+1)$-connected. Since $Y \cap Z \se Z_{(\lan e_n \ran)}$,
$Z$ is $(n-2)$-connected \cite[2.13 (ii)]{kall}. To complete the
proof we have to prove that $Z':=\underline{\u}(F^n)_w$ is
$(n-|w|-2)$-connected. If $w \in Y$, then replacing $Z$ by $Z'$ in the above
and using the induction assumption one sees that $Z'$ is $(n-|w|-2)$-connected.
If $w \notin Y$ then by induction $Y \cap Z'$ is $(n-|w|-2)$-connected
and $Y \cap {Z'}_u$ is $(n-|w|-|u|-2)$-connected for every $u \in Z'\bs Y$
as we proved in the above. Now by \cite[2.13 (i)]{kall} the poset
$Z'$ is $(n-|w|-2)$-connected.
The proof of the last claim is similar to the proof given in
Remark \ref{rem1}(ii).
\end{proof}

\begin{lem}\label{finite2}
Let ${\rm char}(F) \neq {\rm char}(k)$. Then
we have the isomorphism 
$H_i({F^\ast}^p \times G_{n-p}) \simeq H_i(\stab(\sii))$ for all $i$.
\end{lem}
\begin{proof}
Let $M$ be a finite dimensional
$F_1$-vector space, where $F_1:=\{x \in F: \overline{x}=x\}$.
From \cite[Cor. 10.2, Chap. III]{bro} and the fact that
for every group $G$, $H_i(G, k) \simeq {\rm Hom}(H^i(G, k), k)$,
we deduce that
$H_i(M, k)= \begin{cases}
k & \text{if $i=0$}\\
0 & \text{if $i \neq 0$} \end{cases}$.
By a proof similar to the proof of \ref{m-vdk5}, one sees that
$H_i({F^\ast}^p \times G_{n-p}, k) \simeq H_i(\stab(\sii), k)$.
\end{proof}

Applying lemmas \ref{finite1} and \ref{finite2} one sees that
theorem \ref{s-vdk} is true if $F\neq \F_2$ and
${\rm char}(F) \neq {\rm char}(k)$. So we can apply the techniques
that we developed in  sections 3 and 4 to prove the following theorems.

\begin{thm}\label{m6}
Let $F$ be a finite field different from $\F_2$ and
${\rm char}(F) \neq {\rm char}(k)$. Then
\par {\rm (i)} the map $H_l(\inc): H_l(G_{n}) \arr H_l(G_{n+1})$
is surjective for $n \ge l$ and is injective for $n \ge l+1$,
\par {\rm (ii)} if $n-h \ge l$, then the complex
\begin{gather*}
\hspace{-2 cm}
H_l(\rr^h \times G_{n-h}) \overset{\tau_h}{\larr}
H_l(\rr^{h-1} \times G_{n-h+1}) \overset{\tau_{h-1}}{\larr} \cdots \\
\hspace{5 cm}
\overset{\tau_2}{\larr}
H_l(\rr \times G_{n-1}) \overset{\tau_1}{\larr}
H_l( G_{n}) \larr 0
\end{gather*}
is exact, where $\tau_i:=d_{i, l}^1(n)$.
\end{thm}

\begin{thm}
Let ${\rm char}(F)=p$. Then the map
$H_l(\inc): H_l(G_{n}, \z[\frac{1}{p}]) \arr H_l(G_{n+1}, \z[\frac{1}{p}])$
is surjective if $n \ge l+1$ and is injective if $n \ge l+2$.
\end{thm}

\begin{rem}
Let $F$ be a finite field such that ${\rm char}(F) = {\rm char}(k)$.
\par(i) We don't know if a similar results as \ref{m6}
is true or not. There is some information from previous
results, it is true if $n \ge 2l+3$ \cite[Thm. 8.2]{m-vdk}.
\par(ii)
Theorem \ref{m-vdk5} is not true in this case
because otherwise it will be true with every prime field $k$
as a coefficient group
and so it must be true with integral coefficients (see proof of the theorem
\ref{m5}). Hence $\rr^p \times G_{n-p}$ must be isomorphic to the group
$\stab(\sii)$ \cite{cul}, which is not true.
\end{rem}



\begin{thebibliography}{99}

\bibitem{bro} Brown, K. S. Cohomology of groups.  Graduate
Texts in Mathematics, 87. Springer-Verlag, New York, 1994.

\bibitem{ch} Charney, R. A generalization of a theorem of Vogtmann.
J. Pure Appl. Algebra {\bf 44} (1987), 107--125.

\bibitem{cul} Culler, M. Homology equivalent finite groups are isomorphic.
Proc. Amer. Math. Soc. {\bf 72} (1978), no. 1, 218--220.





\bibitem{kar1} Karoubi, M. Relations between algebraic $K$-theory and
Hermitian $K$-theory.
J. Pure Appl. Algebra {\bf 34} (1984), no. 2-3, 259--263.




\bibitem{mac} Mac Lane, S. Homology,
New York; Springer-Verlag, Berlin-G\"ottingen-Heidelberg 1963.

\bibitem{mil} Milnor, J. On the homology of Lie groups made discrete.
Comment. Math. Helv. {\bf 58} (1983), no. 1, 72--85.

\bibitem{m-vdk} Mirzaii, B.; Van der Kallen, W. Homology stability
for unitary groups. Documenta Math.  {\bf 7} (2002) 143--166

\bibitem{nes-sus} Nesterenko Yu. P., Suslin A. A. Homology of the general
linear group over a local ring, and Milnor's $K$-theory.
Math. USSR-Izv. {\bf 34} (1990), no. 1, 121--145.

\bibitem{pan1} Panin I. A. Homological stabilization for the orthogonal
and symplectic groups.
J. Soviet Math. {\bf 52} (1990), no. 3, 3165--3170.

\bibitem{quil6} Quillen D. Characteristic classes of representations.
Algebraic {\it K}-theory, Lecture Notes in Math.,
Vol. {\bf 551}, (1976) 189--216.

\bibitem{sah} Sah, C. Homology of classical Lie groups made discrete.
III. J. Pure Appl. Algebra {\bf 56} (1989), no. 3, 269--312.





\bibitem{kal2} Van der Kallen W. The $K_2$ of rings with many units,
 Ann. Sci. \'Ec. Norm. Sup. (4) {\bf 10} (1977), 473--515.

\bibitem{kall} Van der Kallen W. Homology stability for linear groups.
Invent. Math. {\bf 60} (1980), 269--295.


\bibitem{vog1} Vogtmann, K. Homology stability for ${\rm O}\sb{n, n}$.
Comm. Algebra {\bf 7} (1979), no. 1, 9--38.

\bibitem{vog} Vogtmann, K. Spherical posets and homology stability for
${\rm O}\sb{n, n}$. Topology {\bf 20} (1981), no. 2, 119--132.


\end{thebibliography}
\end{document}